\documentclass[twoside]{article}
\usepackage{amsfonts, amsmath, graphicx, fancyhdr}
\usepackage{float,caption}
\usepackage[top=1.5in, bottom=1.5in, left=1in, right=1in]{geometry}

\newtheorem{lemma}{Lemma}
\newtheorem{remark}{Remark}
\newtheorem{definition}{Definition}
\newtheorem{proposition}{Proposition}
\newtheorem{theorem}{Theorem}

\begin{document}

\title{Generalized solution of a mixed problem for linear hyperbolic system}
\author{Lalla Saadia Chadli, Said Melliani and Aziz Moujahid}
\date{{\small Laboratoire de Mod\'elisation et Calcul (LMC), Facult\'e des Sciences et Techniques} \\
{\small Universit\'e Sultan Moulay Slimane, BP 523, B\'eni Mellal, Morocco} }

\fancyhf{}
\addtolength{\headwidth}{\marginparsep}
\addtolength{\headwidth}{\marginparwidth}
\fancyhead{} 
\fancyhead[RO,LE]{\thepage}
\fancyhead[CE]{Generalized solution of a mixed problem for linear hyperbolic system}
\fancyhead[CO]{L. S. Chadli, S. Melliani and A. Moujahid}
\renewcommand{\headrulewidth}{0.0pt}

\maketitle
 
\begin{abstract}
In the first part of this article, we will prove an existence-uniqueness result for generalized solutions of a mixed problem for linear hyperbolic system in the Colombeau algebra. In the second part, we apply this result to a wave propagation problem in a discontinuous environment.
\end{abstract}

\section{Introduction}

In 1982, Colombeau introduced an algebra $\mathcal{G}$ of generalized functions to deal with the multiplication problem for distributions,see Colombeau \cite{Colombeau1, Colombeau2}. This algebra $\mathcal{G}$ is a differential algebra which contains the space $\mathcal{D}'$ of distributions. Furthermore, nonlinear operations more general than the multiplication make sense in the algebra $\mathcal{G}$. Therefore the algebra $\mathcal{G}$ is a very convenient one to find and study solutions of nonlinear differential equations with singular data and coefficients.

\noindent Consider the mixed problem for the linear hyperbolic system in two variables
\begin{equation}
	\begin{cases}
  	\Bigl( \partial _{t} + \Lambda \left( x,t\right) \partial _{x} \Bigr) U = F(x,t) U + A(x,t) & (x,t) \in (\mathbb{R}_{+}^{*})^{2} \\
    U\left( x,0\right) = U_{0}\left( x\right) & x\in \mathbb{R}_{+} \\
    U_{i}\left( 0,t\right) =\sum\limits_{k=r+1}^{n}\:v _{ik}\left( t\right) 
    U_{k}\left( 0,t\right) +H_{i}\left( t\right) & i=1,\ldots ,r\hspace{0.3in} t\geq 0 \\
    + \text{ Compatibility conditions} &
	\end{cases}   \label{SYS1}
\end{equation}
where $\Lambda $, $F$ and $V$ are $(n\times n)$ matrices whose terms are discontinuous functions. The matrix $\Lambda$ is real and diagonal such that
\[
\Lambda _{1} > \Lambda _{2} > \cdots >\Lambda _{r}>0 > \Lambda _{r+1} > \cdots > \Lambda _{n}
\]
In the case where $\Lambda \in L^{\infty }\left(\mathbb{R}_{+}^{2}\right)$ and $F\in W_{\mbox{loc}}^{-1,\infty }\left(\mathbb{R}_{+}^{2} \right)$, multiplicative products of distributions appear in system~\eqref{SYS1}, and so there is no general way of giving a meaning to system \eqref{SYS1} in the sense of distribution. This hyperbolic system even when it is in the form of a system of conservation laws does not admit any solutions distributions in general see \cite{Hurd}. Our approach is to study \eqref{SYS1} in Colombeau's algebra \cite{Colombeau1, Colombeau2}, and under some hypotheses on $\Lambda $, $F$, $\nu $ and $H$, the system \eqref{SYS1} admits an unique solution in $\mathcal{G}\left( \mathbb{R}_{+}^{2}\right)$. This result completes work already made in the global case by M. Oberguggenberger \cite{Oberguggenberger1}.

\noindent The second part of this article, we will apply this result to the wave propagation problem in a discontinuous environment. the following system
\begin{equation}
	\begin{cases}
		\Bigl( \partial _{t}+c\left( x\right) \partial _{x}\Bigr) u\left(x,t\right) =0 & \left( x,t\right) \in (\mathbb{R}_{+}^{*})^{2} \\
		\Bigl( \partial _{t}-c\left( x\right) \partial _{x}\Bigr) v\left(x,t\right) =0 & \left( x,t\right) \in (\mathbb{R}_{+}^{*})^{2} \\
		u\left( x,0\right) =u_{0}\left( x\right) & x\geq 0 \\
		v\left( x,0\right) =v_{0}\left( x\right) & x\geq 0 \\
		u\left( 0,t\right) =h\left( t\right) v\left( 0,t\right) +b\left( t\right) & t\geq 0 \\
		+ \text{ Compatibility conditions} &
	\end{cases}	\label{SYS2}
\end{equation}
with
\begin{equation*}
	c(x) = \begin{cases}
		c_{R} & \text{if $x > x_0$} \\
		c_{L} & \text{if $0 < x < x_0$}
	\end{cases}
\end{equation*}
$c_{R}$ and $c_{L}$ are real constants, $u_{0}$ and $v_{0}$ are continuous almost everywhere.

\noindent For this problem one can find a classical solution on
\( \left\{0\leq x<x_{0} : t\geq 0\right\} \)
and
\( \left\{ x>x_{0} : t\geq 0\right\} \), \
so imposing a transmission condition in $x=x_{0}$ : the continuity of $u$ and $v$, one will have a classical solution on 
$\left\{ x\geq 0\mbox{ },\mbox{ } t\geq 0 \right\}$.

\noindent Further if $\left(u_{0},v_{0}\right)$ are generalized functions, one can show that the problem \eqref{SYS2} has a unique solution $\left(U,V\right) \in \mathcal{G}\left( \mathbb{R}_{+}^{2}\right) \times \mathcal{G}\left( \mathbb{R}_{+}^{2}\right)$, without having us need of the passage conditions, in the same way one shows that this solution admits an associated distribution that is equal to the classical solution by adjusting.

\section{Existence and uniqueness}

We recall some definitions from the theory of generalized functions which we need in the sequel.

\noindent We define the algebra $\mathcal{G}\left( \mathbb{R}^{m}\right)$ as follows
\[
A_q \left(\mathbb{R}\right) = \Bigl\{ \chi \in \mathcal{D}\left(\mathbb{R}\right) : \int_{\mathbb{R}} \chi (x)\:dx=1 \ \textrm{ and } \ \int_{\mathbb{R}} x^k\:\chi (x)\:dx=0 \quad \textrm{ for } \quad 1\leq k\leq q \Bigr\}
\]
and
\[
A_q \left(\mathbb{R}^m\right) = \Bigl\{ \varphi \left(x_1,\ldots,x_m\right)=\prod^{m}_{j=1}\left\{\chi \left(x_j\right)\right\} \Bigr\}
\]
Let $\mathcal{E}\left[ \mathbb{R}^{m}\right]$ be the set of functions on $\mathcal{A}_1\left( \mathbb{R}^{m}\right)\times \mathcal{C}^{\infty}\left( \mathbb{R}^{m}\right)$ with values in $\mathbb{C}$ witch are $\mathcal{C}^{\infty}$ to seconde variable. Obviously $\mathcal{E}\left[ \mathbb{R}^{m}\right]$ with point wise multiplication is an algebra but $\mathcal{C}^{\infty}\left( \mathbb{R}^{m}\right)$ is not a subalgebra.

\noindent Then given $\varphi \in \mathcal{A}_1\left( \mathbb{R}^{m}\right)$ and $\varepsilon \in \left]0\:,\:1\right[$, we define a function $\varphi_{\varepsilon}$ by
\[
\varphi_{\varepsilon}\left(x\right)=\varepsilon^{-m}\varphi \left(\frac{x}{\varepsilon}\right) \ \ \ \textrm{ for } \ \ x\in \mathbb{R}^{m}
\]
An element of $\mathcal{E}\left[ \mathbb{R}^{m}\right]$ is called "moderate" if for every compact subset $K$ of $\mathbb{R}^{m}$ and every differential operator

\noindent $D=\partial^{k_1}_{x_1},\ldots,\partial^{k_m}_{x_m}$ there is $N\in \mathbb{N}$ such that the following holds
\begin{equation}
	\left\{
\begin{array}{ll}
	\forall\varphi\in \mathcal{A}_N\left( \mathbb{R}^{m}\right), \ \exists C, \exists \eta>0 \ \ \  & \textrm{such that} \\
	\sup\limits_{x\in K}\left|D\:u\left(\varphi_{\varepsilon},x\right)\right|\leq C\varepsilon^{-N} \hspace{0.3in} & \textrm{if} \ 0<\varepsilon<\eta
\end{array}
	\right.
\end{equation}
$\mathcal{E}_M\left[ \mathbb{R}^{m}\right]$ denotes the subset of moderate elements where the index $M$ stands for "\emph{Moderate}". We define an ideal $\mathcal{N}\left[ \mathbb{R}^{m}\right]$ of $\mathcal{E}_M\left[ \mathbb{R}^{m}\right]$ as follows :

\noindent $u\in \mathcal{N}\left[ \mathbb{R}^{m}\right]$ if for every compact subset $K$ of $\mathbb{R}^{m}$ and every differential operator $D$, there is $N\in \mathbb{N}$ such that :
\begin{equation}
	\left\{
\begin{array}{ll}
	\forall q\geq N, \ \forall\varphi\in \mathcal{A}_q\left( \mathbb{R}^{m}\right), \ \exists C, \exists \eta>0 \ \ \  & \textrm{ such that} \\
	\sup\limits_{x\in K}\left|D\:u\left(\varphi_{\varepsilon},x\right)\right|\leq C\varepsilon^{q-N} \hspace{0.3in} & \textrm{if} \ 0<\varepsilon<\eta
\end{array}
	\right.
\end{equation}Finally the algebra $\mathcal{G}\left( \mathbb{R}^{m}\right)$ is defined as the quotient of $\mathcal{E}_M\left[ \mathbb{R}^{m}\right]$ with respect to $\mathcal{N}\left[ \mathbb{R}^{m}\right]$.

In what follows, the elements of $\mathcal{G}\left( \mathbb{R}^{2}\right)$ will be written with capital letters and their representatives in $\mathcal{E}_M\left[ \mathbb{R}^{2}\right]$ with small letters. Furthermore we use the following simplified notations :
	\[
	u\left(\varphi_{\varepsilon},x\right)=u^{\varepsilon}\left(x\right)
\]

\noindent In our work we need a subset of $\mathcal{E}_M\left[ \mathbb{R}_{+}^{2}\right]$ that contains elements $u$ satisfying the following properties :

\begin{itemize}
\item[(a)]  $\exists $ $N\in \mathbb{N}$ such that for all $\varphi \in \mathcal{A}_{N}\left( \mathbb{R}_{+}^{2}\right)$
\[
\exists c>0\hspace{0.3in}\eta >0\mbox{ }:\sup\limits_{y\in \mathbb{R}%
_{+}^{2}}\left| u\left( \varphi _{\varepsilon },y\right) \right| \leq c%
\hspace{0.3in}\textrm{if} \ 0<\varepsilon <\eta
\]
\item[(b)]  For every compact subset $K$ of $\mathbb{R}^{2}_{+}$ , $\exists $ $N\in \mathbb{N}$ such that $\forall \varphi \in \mathcal{A}_{N}\left( \mathbb{R}_{+}^{2}\right)$
\[
\exists c>0\hspace{0.3in}\exists \eta >0\mbox{ }:\sup\limits_{y\in K}\left| u\left( \varphi _{\varepsilon },y\right) \right| \leq N  \log \left( \frac{c}{\varepsilon}  \right)  \hspace{0.3in}\textrm{if} \ 0 <\varepsilon < \eta
\]
\end{itemize}

\begin{definition}
A generalized function $U\in \mathcal{G}\left( \mathbb{R}_{+}^{2}\right) $ admitting a representative $u$ with the property (a) (respectively (b)) is called globally bounded (respectively locally logarithmic growth).
\end{definition}

\begin{definition}
the system \eqref{SYS1} satisfies the compatibility conditions in $\mathcal{G}\left( \mathbb{R}_{+}^{2}\right) $ if there exist $u_{0}^{\varepsilon}$, $\lambda ^{\varepsilon }$, $f^{\varepsilon }$, $h^{\varepsilon }$, $v ^{\varepsilon }$ et $a^{\varepsilon }$ the representatives of $U_{0}$, $\Lambda$, $F$, $H$, $V$ and $A$ that satisfy to the classic conditions compatibility in order to have a $\mathcal{C}^{\infty }$ solution for the classic problem.
\end{definition}

\begin{theorem}
Let $F$, $\Lambda$ and $A$ be $n\times n$ matrices with coefficients in $\mathcal{G}\left( \mathbb{R}_{+}^{2}\right) $, suppose that: there exists $r$ as :
\[
\Lambda _{1} > \Lambda _{2} > \cdots >\Lambda _{r} > 0 > \Lambda _{r+1} > \cdots > \Lambda _{n}
\]
$\Lambda _{i}$ ($i=1,\ldots ,n$) are globally bounded, $\partial _{x}\Lambda _{i}$ and $F_{i}$ are locally logarithmic growth, so for an initial data $U_{0}$ in $\mathcal{G}\left( \mathbb{R}_{+}\right) $, $V _{i}$ an element in $\mathcal{G}\left( \mathbb{R}_{+}\right) $ globally bounded and $H_{i}$ in $\mathcal{G}\left( \mathbb{R}_{+}\right) $, then the problem \ref{SYS1} has an unique solution in $\mathcal{G}\left( \mathbb{R}_{+}^{2}\right) $.
\end{theorem}

\noindent \textbf{Proof} : The proof of the theorem is an adaptation to the demonstration of the theorem 1.2 in \cite{Oberguggenberger1}, therefore one is going to give the big lines rightly.

\noindent Let $\lambda$ a representative of $\Lambda$ in $\mathcal{G} \left( \mathbb{R}_{+}^{+}\right)$ such that
\[
\lambda _{1} > \lambda _{2} > \cdots > \lambda _{r} > 0 > \lambda _{r+1} > \cdots > \lambda _{n}
\]
with $\lambda _{i}$ satisfies the property (a) and $\partial_{x}\lambda _{i}$ satisfies the property (b). \\

\noindent Let $f$ and $a$ are any representatives of $F$ and $A$  in $\mathcal{G}\left( \mathbb{R}_{+}^{2}\right) $ with $f$ satisfies (b).

\noindent $v$, $h$ and $u_{0}$ are any representatives of $V$, $H$ and $U_{0}$ in $\mathcal{G}\left( \mathbb{R}_{+}\right) $ with $v$ satisfies (a).

\noindent so Let's consider the following problem
\begin{equation}
	\begin{cases}
		\Bigl( \partial_{t}+\lambda_{i}^{\varepsilon} (x,t) \partial_{x}\Bigr) u_{i}^{\varepsilon} = \sum\limits_{k=1}^{n} f_{ik}^{\varepsilon} (x,t) u_{k}^{\varepsilon} (x,t) + a_{i}^{\varepsilon} (x,t) & (x,t) \in (\mathbb{R}_{+}^{*})^{2} \\
		u_{i}^{\varepsilon} (x,0) = u_{0_{i}}^{\varepsilon} (x)  & i=1, \ldots, n\quad x\in \mathbb{R}_{+} \\
		u_{i}^{\varepsilon} (0,t) = \sum\limits_{k=r+1}^{n} \nu_{ik}^{\varepsilon} (t) u_{k}^{\varepsilon} (0,t) + h_{i}^{\varepsilon} (t)  & i=1, \ldots, r \quad t\geq 0
	\end{cases} \label{sys5}\tag*{ $\textbf{(I}_{\varepsilon} \textbf{)}$ }
\end{equation}
if we denote $\gamma _{i}^{\varepsilon}$ the corresponding characteristic curve to $\lambda _{i}^{\varepsilon}$ then the problem $\big( \textbf{I}_{\varepsilon} \bigr)$ admits an unique solution $u^{\varepsilon}$, $u_{i}^{\varepsilon}\in \mathcal{C}^{\infty}\left( \mathbb{R}_{+}^{2}\right)$ given by \ \\

\noindent for $i=r+1,...,n$
\begin{eqnarray*}
u_{i}^{\varepsilon }\left( x,t\right) =u_{0_{i}}^{\varepsilon }\left(\gamma ^{\varepsilon}_{i}(x,t,0)
\right) & + & \int_{0}^{t}\Bigl[ \sum\limits_{k=1}^{n}  f_{ik}^{\varepsilon} \Bigl( \gamma _{i}^{\varepsilon }\left( x,t,\tau \right) ,\tau \Bigr)
u_{k}^{\varepsilon }\Bigl( \gamma _{i}^{\varepsilon }\left( x,t,\tau \right)
,\tau \Bigr) \\
  & + & a_{i}^{\varepsilon }\Bigl( \gamma _{i}^{\varepsilon }\left(
x,t,\tau \right) ,\tau \Bigr) \Bigr] d\tau
\end{eqnarray*}
for $i=1,\ldots ,r$%
\begin{eqnarray*}
u_{i}^{\varepsilon }\left( x,t\right)  & = & \sum\limits_{k=r+1}^{n}v
_{ik}^{\varepsilon }\left( t_{0}\right)
\int_{0}^{t_{0}}\sum\limits_{s=1}^{n}\biggl[ f_{ks}^{\varepsilon }\Bigl(
\gamma _{k}^{\varepsilon }\left( 0,t_{0},\tau \right) ,\tau \Bigr)
u_{s}^{\varepsilon }\Bigl( \gamma _{k}^{\varepsilon }\left( 0,t_{0},\tau
\right) ,\tau \Bigr) \biggr] d\tau  \\
& + & \int_{t_{0}}^{t}\sum\limits_{k=1}^{n}\biggl[ f_{ik}^{\varepsilon
}\Bigl( \gamma _{i}^{\varepsilon }\left( x,t,\tau \right) ,\tau \Bigr)
u_{k}^{\varepsilon }\Bigl( \gamma _{i}^{\varepsilon }\left( x,t,\tau \right)
,\tau \Bigr) \biggr] d\tau  \\
& + & \int_{t_{0}}^{t}a_{i}^{\varepsilon }\Bigl( \gamma _{i}^{\varepsilon
}\left( x,t,\tau \right) ,\tau \Bigr) d\tau  \\
& + & \sum\limits_{k=r+1}^{n}v _{ik}^{\varepsilon }\left( t_{0}\right)
\int_{0}^{t_{0}}a_{k}^{\varepsilon } \Bigl( \gamma _{k}^{\varepsilon }\left(
0,t_{0},\tau \right) ,\tau \Bigr) d\tau  \\
& + & \sum\limits_{k=r+1}^{n}v _{ik}^{\varepsilon }\left( t_{0}\right)
u_{0_{k}}^{\varepsilon }\Bigl( \gamma _{k}^{\varepsilon }\left(
0,t_0,0\right) \Bigr) +h_{i}^{\varepsilon }\left( t_{0}\right)
\end{eqnarray*}
where $t_{0}$ is such that the curve $\gamma _{i}$ cuts the axis $%
\left( 0t\right) $ at a point $P_{i}\left( 0,t_{0}\right) $. $%
u_{i}^{\varepsilon }$ is $\mathcal{C}^{\infty }$ function, so it remains to show therefore that $%
u_{i}^{\varepsilon }$ is moderate growth.

\noindent from assumptions, we have
\[
\begin{array}{l}
\exists M>0\quad \mbox{such that : } \displaystyle  \left| \frac{d\gamma _{i}^{\varepsilon
}\left( x,t,\tau \right) }{d\tau }\right| <M\quad \forall (x,t)\in \mathbb{R}_{+}^{2} \quad  \forall i=1,\ldots ,n \\
\exists M_{1}>0\quad \mbox{such that : }\max\limits_{i,j}\left| v
_{i,j}^{\varepsilon }\left( y\right) \right| <M_{1}\quad \forall y\in
\mathbb{R}_{+}
\end{array}
\]
Let $K_{0}$ be a compact in $\mathbb{R}_{+}$, we draw the straight line with a slope $-M$, the determination domain $K_{T}$ of the solution $u_{i}^{\varepsilon }$ does not depend on $\varepsilon $.

\begin{figure}[ht]
\centering
\includegraphics[scale=0.5]{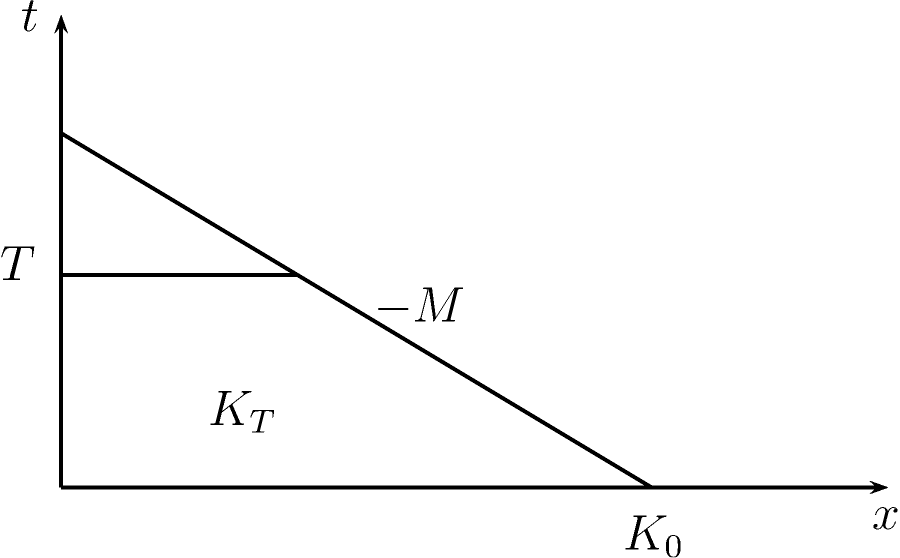}
\caption{}
\label{fig3}
\end{figure}

$\hfill$ $\Box$

\begin{lemma}
Let $u^{\varepsilon }$ a solution of problem ($\mbox{I}_{\varepsilon}$) then $u_{i}^{\varepsilon }$ verified
\begin{eqnarray*}
\sup\limits_{\left( x,t\right) \in K_{T}}\left| u_{i}^{\varepsilon }\left(
x,t\right) \right|  &\leq &M_{2}\left[ \sup\limits_{k}\sup\limits_{\left(
x,t\right) \in K_{T}}\left| a_{k}^{\varepsilon}\left( x,t\right)
\right| .T+\right.  \\
&&\left. \sup\limits_{k}\sup\limits_{x\in K_{0}}\left| u_{0_{k}}^{\varepsilon} \left(x\right) \right| +\sup\limits_{k}\sup\limits_{t\in \left[
0,T\right] }\left| h_{k}^{\varepsilon} \left(t\right) \right| \right] \times
\\
&&\exp \left( nM_{2}\sup\limits_{i,k}\sup\limits_{\left( x,t\right) \in
K_{T}}\left| f_{ik}^{\varepsilon}\left( x,t\right) \right|
.T\right)
\end{eqnarray*}
with
\[
M_{2}=\max \left( nM_{1},1\right)
\]
\end{lemma}
\textbf{Proof} :

\noindent for $i=1,\ldots ,r$, and from the integral equation that verified by $u_{i}^{\varepsilon }$ we have
\begin{eqnarray*}
\sup\limits_{\left( x,t\right) \in K_{T}}\left| u_{i}^{\varepsilon }\left(
x,t\right) \right|  &\leq &M_{2}\left[ T\sup\limits_{\left( x,t\right) \in
K_{T}}\left| a_{k}^{\varepsilon }\left( x,t\right) \right|
+\sup\limits_{k}\sup\limits_{x \in K_{0}}\left|
u_{0_{k}}^{\varepsilon }\left( x\right) \right| +\right.  \\
&&\left. \sup\limits_{k}\sup\limits_{t\in \left[ 0,T\right] }\left|
h_{k}^{\varepsilon }\left( t\right) \right| \right] + \\
&&nM_{2}\int_{0}^{T}\sup\limits_{\left( x,t\right) \in K_{\tau}}\left|
f^{\varepsilon }\left( x,t\right) \right| \sup\limits_{k}\sup\limits_{\left(
x,t \right) \in K_{\tau }}\left| u_{k}^{\varepsilon }\left( x,t
\right) \right| d\tau
\end{eqnarray*}
and the proof is completed by applying the Gronwall's lemma to the function
\[
s\rightarrow \max\limits_{k}\sup\limits_{\left( x,t\right) \in K_{s}}\left|
u_{k}^{\varepsilon }\left( x,t\right) \right|
\]

\noindent For $i=r+1,\ldots ,n$ it is the same way with $t_{0}=0$, $v =0$, $h=0$. \ \\

$\hfill$ $\Box$ \ \\

\noindent the next of the proof of theorem 1, we have \ \\

\noindent \hspace{.5in} $\exists N_1 \in \mathbb{N}$ such that : $\forall \phi \in \mathcal{A}_{N_1}
(\mathbb{R}_+)$
\[
\exists C_1 > 0 \quad \exists \eta > 0: \quad\quad \sup_{(x,t)\in K_T}\left|a^{\varepsilon}(x,t) \right| \leq C_1 \varepsilon^{-N_1} \quad \textrm{if} \ 0<\varepsilon<\eta
\]
\hspace{.5in} $\exists N_2 \in \mathbb{N}$ such that : $\forall \phi \in \mathcal{A}_{N_2}
(\mathbb{R}_+)$
\[
\exists C_2 > 0 \quad \exists \eta > 0: \quad\quad \sup_{x\in K_0}\left|u_0^{\varepsilon}(x) \right| \leq C_2 \varepsilon^{-N_2} \quad \textrm{if} \ 0<\varepsilon<\eta
\]
\hspace{.5in} $\exists N_3 \in \mathbb{N}$ such that : $\forall \phi \in \mathcal{A}_{N_3}
(\mathbb{R}_+)$
\[
\exists C_3 > 0 \quad \exists \eta > 0: \quad\quad \sup_{t\in [0,T]}\left|h^{\varepsilon}(t) \right| \leq C_3 \varepsilon^{-N_3} \quad \textrm{if} \ 0<\varepsilon<\eta
\]
\hspace{.5in} $\exists N_4 \in \mathbb{N}$ such that : $\forall \phi \in \mathcal{A}_{N_4}
(\mathbb{R}^{2}_{+})$
\[
\exists C_4 > 0 \quad \exists \eta > 0: \quad\quad \sup_{(x,t)\in K_T}\left|f^{\varepsilon}(x,t) \right| \leq N_4 \log\left(\frac{C_4}{\varepsilon}\right) \quad \textrm{if} \ 0<\varepsilon<\eta
\]
therefore according to the lemma, we have
\[
\forall \phi \in \mathcal{A}_{N_5}, \ \ \exists C>0,\quad \eta >0: \ \
\sup_{(x,t)\in K_T}\left|u_{i}^ {\varepsilon}(x,t) \right| \leq C_5 \varepsilon^{-N_5} \quad \textrm{if} \ 0<\varepsilon<\eta
\]
with
\[
N_5 = E\left(N_1 +N_2 +N_3 +NTC_4 N_4 \right)+1
\]
for the other derivatives, differentiating the system ($\mbox{I}_{\varepsilon}$) for example with regard to $x$, one gets a system similar to the first. And because $\partial_x \Lambda$ is locally logarithmic growth one gets the same estimation as before, \ldots, then one has
\[
u_{i}^ {\varepsilon} \in \mathcal{E}_{M}(\mathbb{R}^{2}_{+}) \ \ \ \ i=1,\ldots,n
\]
either the existence of the solution for the problem (1) is in  $\mathcal{G}(\mathbb{R}^{2}_{+})$.

\bigskip

\noindent \textbf{Uniqueness}

\medskip

\noindent Let $U$, $V$ two solutions in $\mathcal{G}(\mathbb{R}^{2}_{+})$ of the problem ($\mbox{I}_{\varepsilon}$), with the same initial data and the same boundary values. One must show that so $u^{\varepsilon}$ is a representative of $U$ and
$\mathcal{G}(\mathbb{R}^{2}_{+})$ and if $v^{\varepsilon}$ is a representative of $V$ in
$\mathcal{G}(\mathbb{R}^{2}_{+})$ then $u^{\varepsilon} - v^{\varepsilon} \in \mathcal{N}(\mathbb{R}^{2}_{+})$ see [2].

\noindent indeed :
$u^{\varepsilon} - v^{\varepsilon}$ verifies the same problem that previously and therefore the demonstration is the same. Then one has
\[
u^{\varepsilon} - v^{\varepsilon} = O \left( \varepsilon^q \right)\ \ \ \forall q
\]

$\hfill$ $\Box$

\begin{remark}
To get the solution in the case where $\Lambda \in \textbf{L}^{\infty} (\mathbb{R}^{2}_{+})$, 
$F \in \textbf{W}^{-1,\infty}(\mathbb{R}^{2}_{+})$, one uses the following result. see [4, proposition 2]
\end{remark}

\begin{proposition}
\hspace{3.0in}

{\textbf{a)}} \ Let $\omega\in \textbf{W}^{-1,\infty}_{loc}(\mathbb{R}^{2}_{+})$
then there exist $U\in \mathcal{G}(\mathbb{R}^{2})$ such that: $U$ is associated to
$\omega$ and $U$ is locally logarithmic growth.

\textbf{b)} \ let $\omega\in \textbf{L}^{\infty}(\mathbb{R}^{2})$ then there exist $U\in \mathcal{G}(\mathbb{R}^{2})$ such that: $U$ is associated to $\omega$ and $U$ is
globally bounded, and $\partial^{\alpha}U$ is locally logarithmic growth. \hspace{1.0in}
\(
\alpha = \left(\alpha_1,\alpha_2 \right) \ \ \ \ \mbox{ such that } \ \ \
\left| \alpha \right| = \alpha_1 +\alpha_2  = 1
\)
\end{proposition}

\begin{remark}
For $g \in \textbf{L}^{\infty}(\mathbb{R}_{+})$ one can find $G\in \mathcal{G}(\mathbb{R}_{+})$ such that:
\[ G \approx g \]
and there exist a representative $g^{\varepsilon}$ of $G$ such that $g^{\varepsilon}$
is nil at the neighborhood of $0$ for all $\varepsilon$.
\end{remark}

\medskip

\noindent \textbf{Application}

\medskip

Consider the problem ( \ref{SYS2} )
\begin{equation}
\left\{
\begin{array}{ll}
\Bigl( \partial _{t}+c\left( x\right) \partial _{x}\Bigr) u\left(
x,t\right) =0 & \left( x,t\right) \in (\mathbb{R}_{+}^{*})^{2} \\
\Bigl( \partial _{t}-c\left( x\right) \partial _{x}\Bigr) v\left(
x,t\right) =0 & \left( x,t\right) \in (\mathbb{R}_{+}^{*})^{2} \\
u\left( x,0\right) =u_{0}\left( x\right) & x\geq 0 \\
v\left( x,0\right) =v_{0}\left( x\right) & x\geq 0 \\
u\left( 0,t\right) = v\left( 0,t\right)  &
t\geq 0 \\
+ \textrm{ Compatibility conditions} &
\end{array}
\right.  \nonumber \label{SYS3}
\end{equation}
with
\[
c\left( x\right) =\left\{
\begin{array}{ll}
c_{R} & \mbox{if }x>x_{0} \\
c_{L} & \mbox{if }0<x<x_{0}
\end{array}
\right.
\]

\noindent For the initials data $u_0$, $v_0$ continuous almost everywhere, and nil at neighborhood of $0$.

\noindent the problem (2) admits a classic solution for
\begin{eqnarray*}
	\bigl\{ 0<x<x_0:t\geq0 \bigr\} & \mbox{ and } & \bigl\{ x>x_0 : t\geq0 \bigr\}
\end{eqnarray*}
and while imposing a passage condition on the $x_0$ (continuity of $u$ and $v$ at the point $x_0$ ) then one will have a solution on
\[  \bigl\{ x\geq0 : t\geq0 \bigr\} \]
defined by
\begin{eqnarray*}
	v(x,t) &=& v_0 \left( \gamma_2 (x,t,0) \right) \\
	u(x,t) &=& \left\{ \begin{array}{ll}
	u_0 \left( \gamma_1 (x,t,0) \right) & \mbox{ on (I)} \\
	v(0,t) & \mbox{ on (II)}
	\end{array}
	\right.
\end{eqnarray*}
so one designates by $\Gamma$ the characteristic curve comes from of $\left(0,0\right)$
the part (I) designates the set of $\left(x,t\right)\in \mathbb{R}^{2}_{+}$
below  $\Gamma$. and the part (II) the set the points $\left(x,t\right)$ over $\Gamma$ (see the figure (2)).

$\gamma_1$ the connected curve characteristic corresponding to $c$.

$\gamma_2$ the connected curve characteristic corresponding to $-c$.

\begin{figure}[ht]
\centering
\includegraphics[scale=0.5]{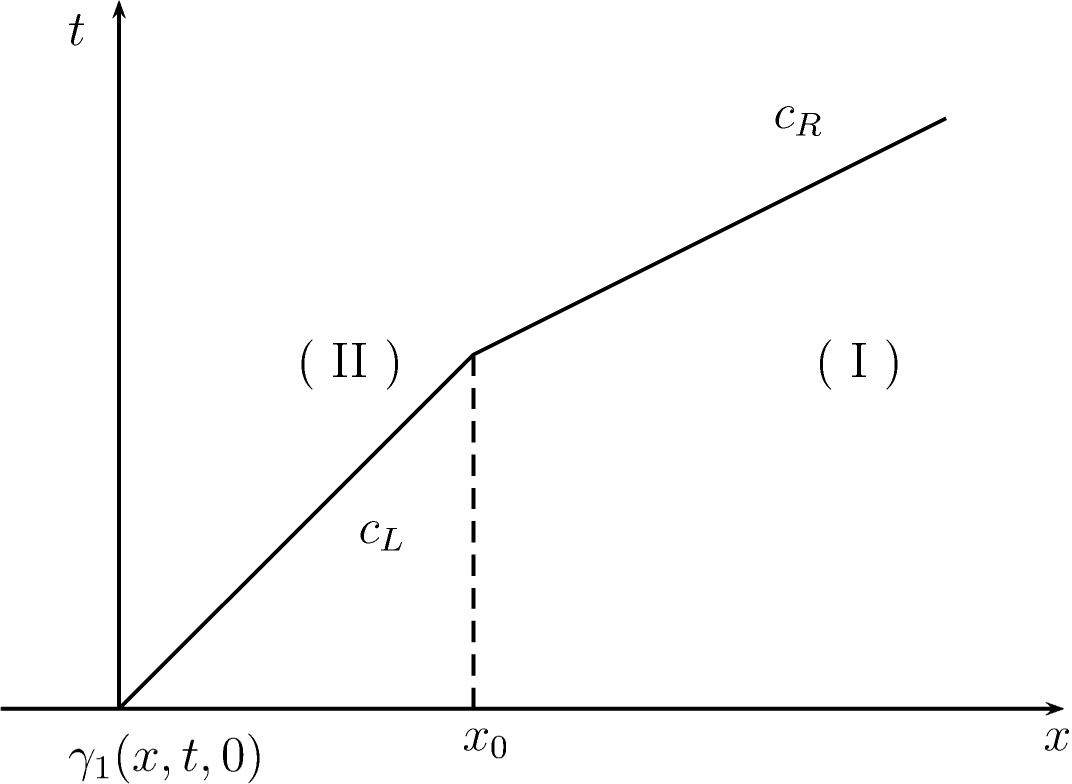}
\caption{}
\label{fig3}
\end{figure}

\begin{proposition}
given $u_0$, $v_0$ two continuous functions nearly everywhere, bounded and nil at the neighborhood of $0$ then the problem (2) admit an unique solution $U$, $V$ in
$\mathcal{G}(\mathbb{R}^{2}_{+})$ besides one has:
\[
U\approx u \hspace{.2in} \mbox{et} \hspace{.2in} V\approx v
\]
with $u$ et $v$ are the distributions solutions of the same problem obtained by
imposing a passage condition.
\end{proposition}

\noindent \textbf{Proof}
$c\in \textbf{L}^{\infty}\left( \mathbb{R}_{+} \right)$, from the proposition (1)
there exists $C\in \mathcal{G}(\mathbb{R}_{+})$ such that
\[ C \approx c \]
$c$ is globally bounded and $\partial_x C$ is locally logarithmic growth.
And so, from the theorem 1, there exists an unique solution $U$, $V$ in $\mathcal{G}(\mathbb{R}^{2}_{+})$ of the problem (2). \ \\

\noindent To show that
\[ U\approx u \]
we suppose that $\left(x,t\right)$ belongs to the region limited by the broken characteristic curve $\Gamma$ comes from the origin and the axis $(ox)$ which we note (region I).

\noindent If $(x,t)$ is over of this curve, the demonstration is identical but with reflection (region II) and for $\left(x,t\right)\in \Gamma$ (the characteristic curve comes from the origin) this set is negligible. \ \\

\noindent let $c^{\varepsilon}$ a representative of $C$ in $\mathcal{G}(\mathbb{R}_{+})$

$u_{0}^{\varepsilon}$ a representative of $U_0$ in $\mathcal{G}(\mathbb{R}_{+})$

$v_{0}^{\varepsilon}$ a representative of $V_0$ in $\mathcal{G}(\mathbb{R}_{+})$

\noindent considering then the following problem
\[
	\left\{
\begin{array}{ll}
	\left(\partial_t + c^{\varepsilon} \partial_x\right)u^{\varepsilon} = 0 &
	\left(x,t\right)\in (\mathbb{R}^{*}_{+})^{2} \\
	\left(\partial_t - c^{\varepsilon} \partial_x\right)v^{\varepsilon} = 0 & \left(x,t\right)\in (\mathbb{R}^{*}_{+})^{2} \\
	u^{\varepsilon} (x,0) = u_{0} ^{\varepsilon}(x) & x\in \mathbb{R}_{+} \\
	v^{\varepsilon} (x,0) = v_{0} ^{\varepsilon}(x) & x\in \mathbb{R}_{+}  \\
    u^{\varepsilon}\left( 0,t\right) = v^{\varepsilon}\left( 0,t\right) & t\in \mathbb{R}_{+}
\end{array}
	\right.
\]
This problem admits an unique solution $u^{\varepsilon}$, $v^{\varepsilon}$ in
$\mathcal{C}^{\infty}(\mathbb{R}^{2}_{+})$. \ \\

\noindent taking

$$\gamma_{1} ^{\varepsilon} = \gamma_1*\phi_{\eta_{\varepsilon}} $$

with $\phi\in \mathcal{D}(\mathbb{R}^{+})$ such that
$$\int_{\mathbb{R}^{+}} \phi(\lambda)d\lambda = 1 \quad\quad supp \ \phi_{\eta_{\varepsilon}} \subset \left]x_{0}-\eta_{\varepsilon},x_{0}+\eta_{\varepsilon} \right[ \quad\quad  \eta_{\varepsilon} = \left| \log \varepsilon \right|^{-1}$$ \ \\

it is evident that for all $(x,t)$ in (region I)
\[
u^{\varepsilon} (x,t) = u_{0} ^{\varepsilon} \bigl(\gamma_{1} ^{\varepsilon}(x,t,0) \bigr)
\]
then to show that $U\approx u$ it is necessary and sufficient to show that :
$\forall\psi\in \mathcal{D}(\mathbb{R}^{2}_{+})$
\[
\lim_{\varepsilon \rightarrow 0} \int_{\mbox{region I}} \Bigl(
u_{0} ^{\varepsilon}\bigl(\gamma_{1} ^{\varepsilon}(x,t,0)\bigr) -
u_{0}\bigl(\gamma_{1}(x,t,0)\bigr) \Bigr)\psi(x,t) dx dt = 0
\]
we have
\begin{eqnarray*}
\int \Bigl( u_{0}^{\:\varepsilon} \bigl( \gamma_{1} ^{\varepsilon} (x,t,0) \bigr) - u_{0} \bigl( \gamma_{1} (x,t,0) \bigr) \Bigr) \psi (x,t) dx\:dt  = \hspace{2.0in} \\
 \int \Bigl( u_{0}^{\varepsilon} \bigl( \gamma_{1} ^{\varepsilon} (x,t,0) \bigr) - u_{0} \bigl( \gamma_{1}^{\:\varepsilon} (x,t,0) \bigr) \Bigr) \psi (x,t) dx\:dt \hspace{1.0in} \\
+  \int \Bigl( u_{0} \bigl( \gamma_{1} ^{\varepsilon} (x,t,0) \bigr) - u_{0} \bigl( \gamma_{1} (x,t,0) \bigr) \Bigr) \psi (x,t) dx\:dt \hspace{1.0in}
\end{eqnarray*}
but
\begin{eqnarray*}
\int \Bigl( u_{0}^{\varepsilon} \bigl( \gamma_{1} ^{\varepsilon} (x,t,0) \bigr) - u_{0} \bigl( \gamma_{1}^{\:\varepsilon} (x,t,0) \bigr) \Bigr) \psi (x,t) dx\:dt \hspace{2.0in} \\
= \int \bigl( u_{0}^{\varepsilon} - u_0 \bigr) \bigl( \gamma_{1} ^{\varepsilon} (x,t,0) \bigr)  \psi (x,t) dx\:dt \hspace{1.0in} \\
\leq  \sup\limits_{x\in \mathbb{R}_+} \left| u_0 \ast \phi_{\varepsilon} - u_0 \right| \left| \int_{\mathbb{R}^{2}_{+}} \psi (x,t) dx\:dt \right| \hspace{1.0in}
\end{eqnarray*}
so
\[
\lim_{\varepsilon \rightarrow 0} \int \Bigl( u_{0}^{\varepsilon} \bigl( \gamma_{1} ^{\varepsilon} (x,t,0) \bigr) - u_{0} \bigl( \gamma_{1}^{\:\varepsilon} (x,t,0) \bigr) \Bigr) \psi (x,t) dx\:dt = 0
\]
to show that
\[
\lim_{\varepsilon \rightarrow 0} \int \Bigl( u_{0} \bigl( \gamma_{1} ^{\varepsilon} (x,t,0) \bigr) - u_{0} \bigl( \gamma_{1} (x,t,0) \bigr) \Bigr) \psi (x,t) dx\:dt = 0
\]
it is sufficient to show that
\[
\lim_{\varepsilon \rightarrow 0} \Bigl( \gamma_{1} ^{\varepsilon} (x,t,0) - \gamma_{1} (x,t,0) \Bigr) = 0
\]
or $c$ is globally bounded, then
\[
\exists M >0 \quad\quad \sup\limits_{x\in \mathbb{R}_+} \left| c^{\varepsilon}(x) \right| < M
\]
so we can to surround the curve $\gamma_{1} ^{\varepsilon}$ between two broken curves, (see the figure 3 ).


\noindent and taking the intersection of these two curves with the axis $(0x)$, it gives us two points
\begin{eqnarray*}
	x_1 & = & c_L \Bigl( - \frac{2 \eta_{\varepsilon}}{M} - \frac{x_0 + \eta_{\varepsilon} - x}{c_R} - t \Bigr) - \eta_{\varepsilon} + x_0 \\
	x_2 & = & - c_L \Bigl( - \frac{2 \eta_{\varepsilon}}{M} + \frac{x_0 + \eta_{\varepsilon} - x}{c_R} + t \Bigr) - \eta_{\varepsilon} + x_0
\end{eqnarray*}
such that
\[
x_1 \leq \gamma_{1} ^{\varepsilon} (x,t,0) \leq x_2
\]
\begin{figure}[ht]
\centering
\includegraphics[scale=0.6]{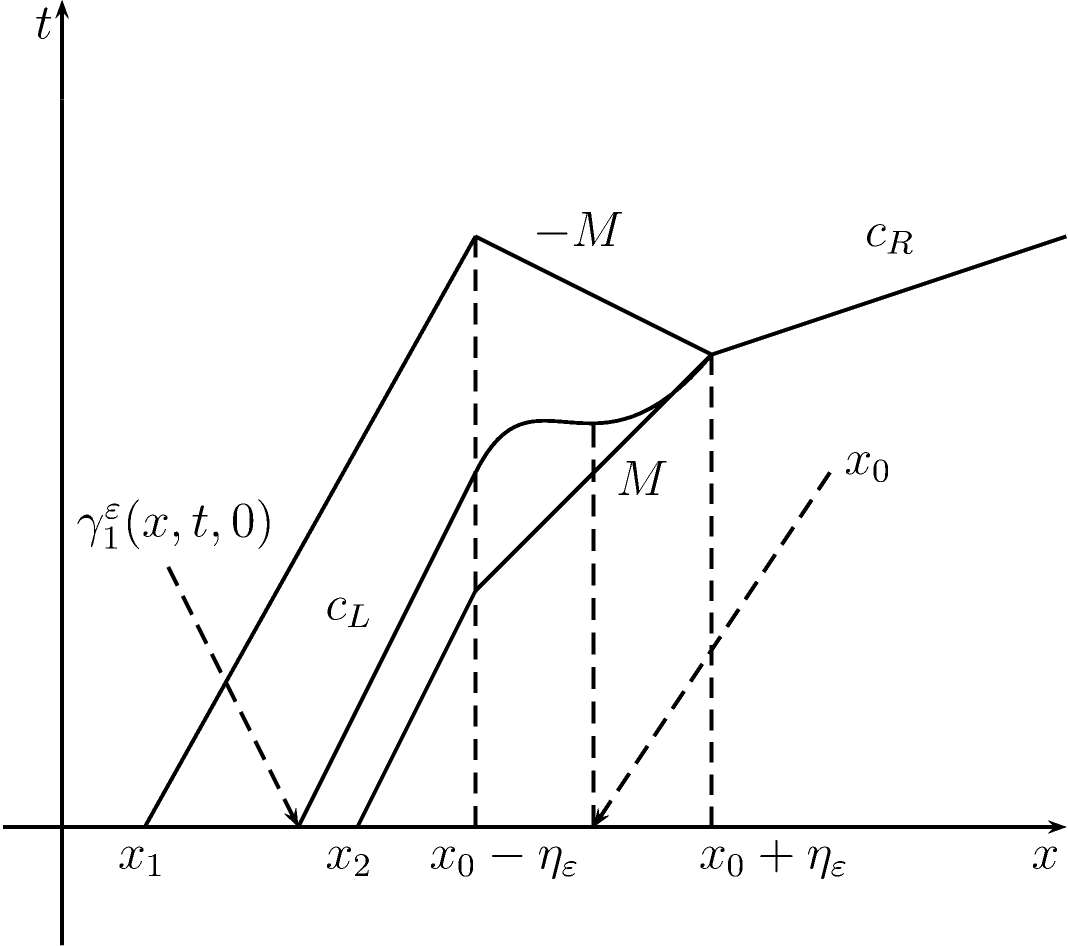}
\caption{}
\label{fig3}
\end{figure}


hence
\begin{eqnarray*}
	\lim_{\varepsilon \rightarrow 0} \gamma_{1} ^{\varepsilon} (x,t,0) &=& - c_L t + \frac{c_L}{c_R} \bigl( x - x_0 \bigr) + x_0 \\
	&=&  \gamma_{1} (x,t,0)
\end{eqnarray*}
then
\[
U \approx u
\]
for $v$, the demonstration is the same.

 $\hfill$ $\Box$

\end{document}